\documentclass[11pt,leqno]{article}

\usepackage{amsfonts,latexsym, amsmath, amssymb, amsthm, hyperref}

\usepackage{gfsartemisia}
\usepackage[T1]{fontenc}

\usepackage{color}
\usepackage{fancyhdr} \usepackage{lastpage}

\newtheorem{theorem}{Theorem}[section]
\newtheorem{lemma}[theorem]{Lemma}

\newtheorem{proposition}[theorem]{Proposition}
\newtheorem{corollary}[theorem]{Corollary}

\theoremstyle{definition}

\theoremstyle{remark}
\newtheorem{remark}[theorem]{Remark}
\newtheorem{note}[theorem]{Note}
\newtheorem*{note*}{Note}

\makeatletter

\numberwithin{equation}{section}

\makeatother

\DeclareMathOperator{\signum}{sgn}

\newcommand{\sgn}[1]{\signum(#1)}
\newcommand{\prend}{$\hfill \quad \Box$}

\newcommand\blfootnote[1]{%
  \begingroup
  \renewcommand\thefootnote{}\footnote{#1}%
  \addtocounter{footnote}{-1}%
  \endgroup
}

\begin{document}

\small

\title{On Dvoretzky's theorem for subspaces of $L_p$
%\footnote{Created: August 1, 2015, Columbia MO.}
}

\author{Grigoris Paouris\thanks{Supported by the NSF CAREER-1151711 grant;} \and Petros Valettas\thanks{Supported by the NSF grant DMS-1612936.} }

%\date{}

\maketitle

\begin{abstract}\footnotesize
We prove that for any $2<p<\infty$ and for every $n$-dimensional subspace $X$ of $L_p$, represented on $\mathbb R^n$, 
whose unit ball $B_X$ is in Lewis' position one has 
the following two-level Gaussian concentration inequality:
\[ \mathbb P\left( \big| \|Z\| - \mathbb E\|Z\| \big| > \varepsilon \mathbb E\|Z\| \right) \leq 
C \exp \left (- c \min \left\{ \alpha_p \varepsilon^2 n, (\varepsilon n)^{2/p} \right\} \right), \quad 0<\varepsilon<1 , \] 
where $Z$ is a standard $n$-dimensional Gaussian vectors, $\alpha_p>0$ is a constant depending only on $p$ and $C,c>0$ are absolute constants. 
As a consequence we show optimal lower bound for the dimension of almost spherical sections for these spaces.
In particular, for any $2<p<\infty$ and every $n$-dimensional subspace $X$ of $L_p$, the Euclidean
space $\ell_2^k$ can be $(1+\varepsilon)$-embedded into $X$ with $k\geq c_p \min\{ \varepsilon^2 n , (\varepsilon n)^{2/p}\}$, 
where $c_p>0$ is a constant depending only on $p$. This improves upon the previously known estimate due to Figiel, Lindenstrauss
and V. Milman.
\end{abstract}


%%%%%%%%%%%%%%%%%%%%%%%%%%%%%%%%%%%%%%%%%%%%%%
\section{Introduction}
%%%%%%%%%%%%%%%%%%%%%%%%%%%%%%%%%%%%%%%%%%%%%%

In the present note we study the classical result of Dvoretzky \cite{Dvo} on almost spherical sections of normed spaces in the case
of subspaces of $L_p$. Grothendieck in \cite{Gro} motivated by the well
known Dvoretzky-Rogers lemma from \cite{DR} to ask if every finite-dimensional normed space 
has lower dimensional subspaces which are almost Euclidean and their dimension grows with respect to the 
dimension of the ambient space. Dvoretzky in \cite{Dvo} gave an affirmative answer in the above question by proving that
for any positive integer $k$ and every $\varepsilon \in (0,1)$ there exists $N=N(k,\varepsilon)$ with the following property:
For every $n\geq N$ and any $n$-dimensional normed space $X$ there exists $k$-dimensional subspace $E$ which 
$(1+\varepsilon)$-isomorphic to the Euclidean space $\ell_2^k$.
In modern functional analytic language 
this means that every 
infinite-dimensional Banach space contains $\ell_2^n$'s uniformly. Dvoretzky's proof in \cite[Theorem 1]{Dvo}
provides the quantitative estimate $N(k,\varepsilon) \geq \exp(c \varepsilon^{-2} k^2 \log^2 k )$ (see \cite{Sza} for a related discussion), for some absolute constant $c>0$\footnote{Here and elsewhere in this paper $c$ and $C$ denote positive absolute constants.} 
However, the aforementioned estimate is not optimal. The optimal dependence with respect to the dimension was proved later by V. Milman 
in his groundbreaking work \cite{Mil} where he obtained $N(k,\varepsilon) \geq \exp(c k \varepsilon^{-2} \log\frac{1}{\varepsilon} )$ (an alternative approach which yields the same estimate was presented by Szankowski in \cite{Sza}).
Equivalently, this states that for any $\varepsilon \in (0,1)$ there exists a function $c(\varepsilon)>0$ with
the following property: for every $n$-dimensional normed space $X$ there exists $k \geq c(\varepsilon) \log n$ 
and a linear map $T:\ell_2^k \to X$ with $\|x\|_2 \leq \|Tx\|_X \leq (1+\varepsilon) \|x\|_2$ for all $x \in \ell_2^k$. In that case we say
that $\ell_2^k$ can be $(1+\varepsilon)$-embedded into $X$ or that $X$ has a $k$-dimensional subspace which is $(1+\varepsilon)$-Euclidean 
and we write $\ell_2^k \stackrel{1+\varepsilon} \hookrightarrow X$. 

The example of $X= \ell_\infty^n$ shows that this result is best possible with respect to $n$ (see \cite{Mil} or \cite[Proposition 3.2]{FLM} for the details). 
The approach of \cite{Mil} is probabilistic in nature and provides that the vast majority of subspaces (in terms of the Haar probability measure on
the Grassmannian manifold $G_{n,k}$) are $(1+\varepsilon)$-spherical, as long as $k\leq c(\varepsilon) k(X)$, where 
$k(X)$ is the {\it critical dimension} of $X$ (see below for the definition). Nowadays this is customary addressed as the randomized Dvoretzky theorem or random
version of Dvoretzky's theorem. V. Milman in this work revealed the significance of the concentration of measure as a basic tool for the 
understanding of the high-dimensional structures.  That was the starting point for many applications of the concentration of measure method
in high-dimensional phenomena. Since then, this tool has found numerous applications in various fields such as quantum information \cite{ASW}, 
combinatorics \cite{BLM}, random matrices \cite{Ver}, compressed sensing \cite{FR}, theoretical computer science \cite{Ma}, 
geometry of high-dimensional probability measures \cite{EK} and more.

Another remarkable fact of V. Milman's approach is that the critical quantity $k(X)$ can
be described in terms of the global parameters of the space. In particular, $k(X) \simeq \mathbb E\|Z\|_X^2 / b^2(X)$ where $Z$ is a standard
Gaussian random vector in $X$ and $b(X)=\max_{\|\theta \|_2=1}\|\theta\|_X$. Then, one can find a good position of the unit ball of
$X$ for which $k(X)$ is large enough with respect to $n$ (see \cite{MS} for further details). It has been proved in \cite{MS2} that this 
formulation is optimal with respect to the dimension $k(X)$ in the sense that the $k$-dimensional subspaces which 
are $4$-Euclidean with probability greater than $\frac{n}{n+k}$ cannot exceed $Ck(X)$ (see \cite{HW} for a recent development 
on this fact).

The proof of \cite{Mil} gave the estimate
$c(\varepsilon) \geq c \varepsilon^2 /\log \frac{1}{\varepsilon}$ and this was improved to $c(\varepsilon) \geq c\varepsilon^2$ 
by Gordon in \cite{Go} and later, adopting the methods of V. Milman, by Schechtman in \cite{Sch1}. This dependence is known to 
be optimal in the setting of the randomized Dvoretzky theorem; see \cite{Sch4}). The works of Schechtman in \cite{Sch3} and Tikhomirov in \cite{Tik} established that the dependence 
on $\varepsilon$ in the randomized Dvoretzky for $\ell_\infty^n$ is of the order $\varepsilon/ \log \frac{1}{\varepsilon}$ and this is 
best possible. Optimal bounds on $c(\varepsilon)$ in the randomized Dvoretzky for $\ell_p^n, \; 1 \leq p \leq \infty$ have recently been studied in \cite{PVZ}.

As far as the dependence on $\varepsilon$ in the ``existential version" of Dvoretzky's theorem is concerned, 
Schechtman proved in \cite{Sch2} that one can always $(1+\varepsilon)$-embed $\ell_2^k$ in any
$n$-dimensional normed space $X$ with $k \geq c\varepsilon \log n / (\log\frac{1}{\varepsilon})^2 $. Tikhomirov in \cite{Tik} 
proved that for 1-symmetric space $X$ we may have $k\geq c \log n/ \log \frac{1}{\varepsilon}$ and this was subsequently extended by
Fresen in \cite{Fres} for permutation invariant spaces with bounded basis constant. For more detailed information on the subject, explicit statements and historical remarks the reader is referred to the recent monograph \cite{AGM}.

The purpose of this note is to study the dependence on $\varepsilon$ and dimension in Dvoretzky's theorem for finite-dimensional 
subspaces of $L_q, \, 2< q <\infty$. The case of subspaces of $L_q, \; 1\leq q < \infty$ have been previously studied in the 
classical paper \cite{FLM} by Figiel, Lindenstrauss and V. Milman.

%%%%%%%%%%%%%%%%%
The approach in \cite{FLM} is based on V. Milman's
asymptotic formula and the fact that the $L_p$ spaces enjoy the {\it cotype} property.
Let us recall that for $2\leq q<\infty$ the $q$-cotype constant of a normed space $X$ in $n$ vectors, denoted
by $C_q(X,n)$, is defined as the smallest constant $C>0$ which satisfies:
\begin{align*}
 \left( \sum_{i=1}^n \|z_i\|_X^q \right)^{1/q} \leq C \mathbb E\left \| \sum_{i=1}^n \varepsilon_i z_i\right \|_X,
\end{align*} for any $n$ vectors $z_1,\ldots,z_n\in X$. Then, the $q$-cotype constant of $X$ is defined as 
$C_q(X):=\sup_n C_q(X,n)$. 
% It is said that $X$ is of cotype $q$ (or it has cotype $q$) if $C_q(X) <\infty$.
Following the terminology of G. Pisier, the notion of cotype is a {\it super-property}, that is it depends only on the finite
dimensional subspaces of the space. It is also isomorphic invariant and the spaces $L_p, \; 1\leq p<\infty$ are of cotype $q=\max\{2,p\}$
with $C_q(L_p) =O(q^{1/2})$ (see \cite{AK} for a proof). Therefore, for any finite-dimensional $X$ of $L_q, \; 2<q<\infty$ we have
$C_q(X)\leq C\sqrt{q}$. The authors in \cite{FLM}, using the the classical Dvoretzky-Rogers lemma, show that any $n$-dimensional 
normed space $X$ of cotype $q$ whose unit ball is in John's position (see e.g. \cite{FLM} for the related definition) 
satisfies $k(X)\geq cC_q^{-2}(X) n^{2/q}$. It follows that if $X$ is an $n$-dimensional
subspace of $L_q, \; 2<q<\infty$ whose unit ball is in John's position then $k(X) \geq c q^{-1}n^{2/q}$ and the 
standard concentration techniques yield $(1+\varepsilon)$-spherical sections of $X$ of 
dimension $k\geq cq^{-1}\varepsilon^2 n^{2/q}$ (see \cite{FLM} for the details).
Moreover, the same argument provides $k(X)\geq cn$ for any 
$n$-dimensional subspace $X$ of $L_q$ with $1\leq q < 2$ in John's position, 
and thus $\ell_2^k$ can be $(1+\varepsilon)$-embedded into $X$ with
$k\geq c \varepsilon^2 n$ which is best possible. In the present note we show that for the range $2<q<\infty$ the estimate can
be considerably improved. Our result reads as follows.
%%%%%%%%%%%%%%%%

\begin{theorem} \label{thm:main-1-1}
 For any $2<p<\infty$ there exists a constant $c(p)>0$ with the following property: 
 for any $n$-dimensional subspace $X$ of $L_p$ and for any $\varepsilon \in (0,1)$
 there exists $k\geq c(p) \min\{ \varepsilon^2 n, (\varepsilon n)^{2/p}\}$ so that 
 $\ell_2^k$ can be $(1+\varepsilon)$-embedded into $X$.
\end{theorem}

Our approach is different and depends on a Gaussian functional analytic inequality rather than the spherical 
isoperimetric inequality that is used in the classical framework. Thus, our proof still depends on random methods,
but the main tool is a variant of an inequality due to Pisier from \cite{Pis}.

To prove the above theorem, we have to bypass Milman's asymptotic formula, 
which involves the Lipschitz constant of the norm. As several examples show this parameter is insufficient to describe 
efficiently phenomena in the almost isometric scale. Our argument outclasses the latter one since it takes 
into account the order of magnitude of the length of the gradient of the norm.
The idea of sufficiently estimating averages of the Euclidean norm of the gradient of a function in order to get
sharp concentration results seems to be only recently applied and was also successfully exploited in \cite{PVZ}.
Moreover, the selection of the position of the unit ball of the space is 
different. Instead of using John's position we employ Lewis' position (see Section 2 for details) for the unit 
ball of finite-dimensional subspaces of $L_p$. This permits us to express the norm in an integral form,
with respect to some isotropic measure on the sphere, and therefore to use the aforementioned inequality. In fact we derive
Theorem \ref{thm:main-1-1} from the randomized Dvoretzky theorem for those spaces in Lewis' position.
For this end we prove that the norm of the underlying subspace in this position
exhibits two-level Gaussian concentration and minimal fluctuations.

\begin{theorem} \label{thm:main-1-2}
Let $2<p<\infty$ and let $X$ be an $n$-dimensional subspace of $L_p$ represented on $\mathbb R^n$ 
whose unit ball $B_X$ is in Lewis' position. Then, 
\[ \mathbb P\left( \big| \|Z\| - \mathbb E\|Z\| \big| > \varepsilon \mathbb E\|Z\| \right) \leq 
C \exp \left (- c \min \left\{ \alpha_p \varepsilon^2 n, (\varepsilon n)^{2/p} \right\} \right), \quad 0<\varepsilon<1 . \] 
In particular, we have
\[ {\rm Var} \|Z\| \leq C_p n^{2/p-1} , \]
where $\alpha_p,C_p>0$ are constants depending only on $p$ and $Z$ is a standard $n$-dimensional Gaussian vector.
\end{theorem}

It is worth mentioning that the
Gaussian concentration and the variance estimate obtained for these spaces is best possible (up to constants of $p$) 
as the example of $\ell_p$ norms shows (see \cite{PVZ} for the exact formulation). 
Consequently, the random version of Dvoretzky's theorem we prove for this position (or for this type of norms) is sharp in the sense that 
in the case of $\ell_p^n$ spaces the corresponding critical dimension is optimal (see \cite{PVZ}). In other words the $\ell_p^n$ space occurs
as the approximately extremal structure in this study or is the worst subspace of $L_p$ with respect to the local almost Euclidean structure.

The novelty of the work is not only observed in the techniques used but also in the content of the results. For our analysis is crucial
the perspective 
% of replacing John's position by Lewis' for the unit ball of the ...
of differently selecting the position of the unit ball of the underlying space and this is reflected in the improved estimates we obtain.
To the best of our knowledge the concentration estimates we derive in Theorem \ref{thm:main-1-2} are new and it is also clear
that the dimension $k(n, p, \varepsilon) \simeq_p \min\{ \varepsilon^2 n , (\varepsilon n)^{2/p}\}$, that one can find almost Euclidean subspaces, 
is always better than the previously known 
$\varepsilon^2 n^{2/p}$ due to Figiel, Lindestrauss and V. Milman. In addition, the improved estimate for $k( n, p, \varepsilon)$ yields
``new dimensions" of almost Euclidean sections in the following sense: The previous setting was only permitting almost isometric 
embeddings of distortion $1+\varepsilon$ with $\varepsilon \gg n^{-1/p}$ in order to achieve non-trivial dimensions. Now this phenomenon
admits an improvement and one can find $(1+\varepsilon)$-linear embeddings with $\varepsilon \gg n^{-1/2}$. It is worth mentioning 
that the dimension $k(n,p,\varepsilon)$ that one finds almost Euclidean sections for these spaces is given
implicitly as function of $\varepsilon$ and $n$ rather than as function of separated variables as V. Milman's formula suggests. This phenomenon 
had not been observed prior to this work and \cite{PVZ}.

The rest of the paper is organized as follows: In Section 2 we introduce the notation, some background material on isotropic
measures on the $(n-1)$-dimensional Euclidean sphere and
finally we give the proof of the aforementioned Gaussian inequality. In Section 3 we prove concentration results for 
the family of the
$L_q$-bodies associated with an isotropic measure $\mu$ on the $(n-1)$-dimensional Euclidean sphere. 
In Section 4 we provide the proof of our main result. Finally, in Section 5 we conclude with some further remarks.

%%%%%%%%%%%%%%%%%%%%%%%%%%%%%%%%%%%%%%%%%%%%%%
\section{Background material and auxiliary results}
%%%%%%%%%%%%%%%%%%%%%%%%%%%%%%%%%%%%%%%%%%%%%%

We work in $\mathbb R^n$ equipped with the standard Euclidean structure $\langle \cdot, \cdot \rangle$. The $(n-1)$-dimensional
Euclidean sphere is defined as $S^{n-1}:=\{x\in \mathbb R^n : \langle x,x\rangle=1\}$. The $\ell_p$ norm
is defined as $\|x\|_p: = (\sum_{i=1}^n |x_i|^p )^{1/p}$ for $x=(x_1, \ldots, x_n) \in \mathbb R^n$. We set 
$\ell_p^n=(\mathbb R^n, \|\cdot\|_p)$ and let $B_p^n$ its unit ball. More generally, for any centrally symmetric
convex body $K$ on $\mathbb R^n$ we write $\|\cdot\|_K$ for the norm induced by $K$. The $n$-dimensional Lebesgue
measure (volume) of a body $A$ is denoted by $|A|$. The space $L_p(\Omega,{ \cal E}, \mu), \; 1\leq p<\infty$ consists 
of all ${\cal E}$-measurable functions $f:\Omega \to \mathbb R$ so that $\int_\Omega |f|^p \, d\mu<\infty$, equipped with the 
norm $\|f\|_{L_p(\mu)} := (\int_\Omega |f|^p \,d \mu)^{1/p}$.

\noindent The $n$-dimensional (standard) Gaussian measure is denoted by $\gamma_n$ and its density is
\begin{align*}
d\gamma_n(x) := (2\pi)^{-n/2} e^{-\|x\|_2^2/2} dx.
\end{align*}
More generally, let $d\gamma_{n,\sigma}(x) := (2\pi \sigma^2)^{-n/2} e^{-\|x\|_2^2/(2\sigma^2)} dx$ for $\sigma>0$. 
Random vectors, usually distributed according to 
$\gamma_n$, are denoted by $Z,W \ldots$ while the random variables by $g_i,\xi, \ldots$. 
The notation $\mathbb E(\cdot)$ is used for the expectation. 
The moments with respect to $\gamma_n$ of norms whose unit ball is the body $K$ are denoted by
\begin{align*}
I_r(\gamma_n,K) := \left( \mathbb E\|Z\|_K^r \right)^{1/r} = \left( \int_{\mathbb R^n} \|z\|_K^r \, d\gamma_n(z) \right)^{1/r}
\end{align*} and more generally, for an arbitrary probability measure $\nu$ as $I_r(\nu, K)$.
Recall the $p$th moment $\sigma_p$ of a standard Gaussian 
random variable $g$
\begin{align} \label{eq: 2.3}
\sigma_p^p := \mathbb E|g|^p =\frac{2^{p/2}}{\sqrt{\pi}} \Gamma \left(\frac{p+1}{2} \right) 
\sim \sqrt{2/e} \left( \frac{p+1}{e}\right)^{p/2}, \quad p\to\infty,
\end{align} where $f\sim h$ means $f(t)/h(t)\to 1$ as $t \to \infty$. We write $f\lesssim h$ when there exists absolute constant
$C>0$ such that $f \leq Ch$. We write $f\simeq h$ if $f\lesssim h$ and $h \lesssim f$, whereas the notation $f\lesssim_p h$ means
that the involved constant depends only on $p$. The letters $C,c, C_1, c_0, \ldots$ are frequently used throughout the text in order
to denote absolute constants which may differ from line to line. 

The random version of Dvoretzky's theorem due to V. Milman from \cite{Mil} (for the optimal dependence on $\varepsilon$
see \cite{Go} and \cite{Sch1}) reads as follows.

\begin{theorem}
Let $X=(\mathbb R^n, \|\cdot\|)$ be a normed space. Define the critical dimension of $X$ as 
the quantity
\begin{align*}
k(X): = \frac{\mathbb E\|Z\|^2}{b^2(X)}, \quad Z\sim N({\bf 0},I_n)
\end{align*} where $b(X):=\max_{\theta \in S^{n-1}}\|\theta\|$. Then, for every $\varepsilon\in (0,1)$ and 
for any $k\leq c\varepsilon^2 k(X)$ the random (with respect to the Haar measure on the Grassmannian $G_{n,k}$) $k$-dimensional subspace $F$ 
of $X$ is $(1+\varepsilon)$-spherical, i.e.
\begin{align*}
\frac{1-\varepsilon}{M} B_F \subseteq B_X \cap F \subseteq \frac{1+\varepsilon}{M} B_F,
\end{align*} with probability greater than $1-e^{-ck(X)}$, where $M=M(X)=\int_{S^{n-1}} \|\theta\| \, d\sigma(\theta)$ and $\sigma$ is the uniform 
probability measure on $S^{n-1}$.
\end{theorem}

%%%%%%%%%%%%%%%%%%%%%%%%%%%%%%%%%%%%%%%%%%%%%%%
\subsection{Logarithmic Sobolev inequality}
%%%%%%%%%%%%%%%%%%%%%%%%%%%%%%%%%%%%%%%%%%%%%%%

Let $\nu$ be a Borel probability measure on $\mathbb R^n$ which satisfies a log-Sobolev 
inequality with constant $\rho>0$
\begin{align*}
{\rm Ent}_\nu(f^2):= \int f^2 \log f^2 \, d\nu - \int f^2\, d\nu \log \left( \int f^2 \, d\nu \right) \leq \frac{2}{\rho} \int_{\mathbb R^n} \| \nabla f \|_2^2\, d\nu,
\end{align*} for all smooth (or locally Lipschitz) functions $f:\mathbb R^n \to \mathbb R$. 
The $n$-dimensional Gaussian measure satisfies the log-Sobolev inequality with $\rho=1$ (see \cite{Led}). The next 
lemma is essentially from \cite{SV} (see also \cite{PVZ}). We provide a sketch of proof for reader's convenience.

\begin{lemma}\label{lem:conc-log-sob}
Let $\nu$ be a Borel probability measure on $\mathbb R^n$ which satisfies a log-Sobolev inequality with constant $\rho$. Then, for any smooth
function $f:\mathbb R^n \to \mathbb R$ we have
\begin{align} \label{eq: 2.7}
\|f\|_{L_q(\nu)}^2 - \|f\|_{L_p(\nu)}^2 \leq \frac{1}{\rho} \int_p^q \big\| \, \|\nabla f\|_2 \, \big \|_{L_s(\nu)}^2 \, ds,
\end{align} for all $2\leq p\leq q$. Moreover, if $f$ is Lipschitz continuous, then we have
\begin{align*}
\|f\|_{L_q(\nu)}^2 - \|f\|_{L_p(\nu)}^2 \leq \frac{\|f\|_{\rm Lip}^2 }{\rho}(q-p) .
\end{align*}
In particular, we obtain
\begin{align*}
\frac{\|f\|_{L_q(\nu)}}{\|f\|_{L_2(\nu)}}\leq \sqrt{1+\frac{q-2}{\rho k(f)}},
\end{align*} for $q\geq 2$, where $k(f):=\|f\|^2_{L_2(\nu)}/\|f\|_{\rm Lip}^2$.
\end{lemma}

\noindent {\it Sketch of Proof.} For $p\geq 2$ we define $I(p):= \|f\|_{L_p}$. Differentiation with respect to $p$ yields
\begin{align*}
\frac{dI}{dp} = \frac{{\rm Ent}_\nu (|f|^p)}{p^2 I(p)^{p-1}}.
\end{align*} Applying the log-Sobolev inequality for $g=|f|^{p/2}$ we obtain
\begin{align*}
\frac{dI}{dp} \leq \frac{1}{2 \rho I(p)^{p-1}} \int_{\mathbb R^n} |f|^{p-2} \| \nabla f \|_2^2 \, d\nu \leq 
\frac{1}{2\rho I(p)^{p-1}} I(p)^{p-2} \big\| \|\nabla f\|_2 \big \|_{L_p(\nu)}^2, 
\end{align*} by H\"{o}lder's inequality. This shows that $(I(p)^2)' \leq \frac{1}{\rho} \big\| \|\nabla f\|_2 \big \|_{L_p(\nu)}^2$, 
thus integration over the interval $[p,q]$ proves \eqref{eq: 2.7}. \prend

%%%%%%%%%%%%%%%%%%%%%%%%%%%%%%%%%%%%%%%%%%%%%%%%%%%
\subsection{Lewis' position}
%%%%%%%%%%%%%%%%%%%%%%%%%%%%%%%%%%%%%%%%%%%%%%%%%%%

Given any finite Borel measure $\mu$ on $S^{n-1}$ (which is not supported in any hyperplane) and any $1\leq p<\infty$ 
we can equip $\mathbb R^n$ with the norm 
\begin{align*}
\|x\|_{\mu,p}: = \left( \int_{S^{n-1}} |\langle x,\theta \rangle|^p \, d\mu(\theta) \right)^{1/p}.
\end{align*} It's clear that the space $X=(\mathbb R^n, \|\cdot\|_{\mu,p} )$ can be naturally embedded into $L_p(S^{n-1},\mu)$ via the
linear isometry $U: X \to L_p(S^{n-1},\mu)$ with $Ux: = \langle x,\cdot \rangle$.

\smallskip

Lewis' fundamental result from \cite{Lew}, states that the previous situation can always be realized for finite-dimensional 
subspaces of $L_p(\nu)$ after a suitable change of the density $\nu$
(see also \cite{SZ} for an alternative proof which extends to the whole range $0<p<\infty$ and arises as a solution of an optimization problem). 
The formulation we use here follows the exposition from \cite{LYZ}.

\begin{theorem} [Lewis] \label{thm:Lew_LYZ}
Let $1\leq p<\infty$ and let $X$ be an $n$-dimensional subspace of $L_p$. Then, there exists an even Borel measure $\mu$ on $S^{n-1}$ 
which satisfies
\begin{align} \label{eq:iso}
\|x\|_2 ^2 =\int_{S^{n-1}} |\langle x, \theta\rangle|^2 d\mu(\theta),
\end{align} for all $x\in \mathbb R^n$ and the normed space $(\mathbb R^n , \|\cdot\|_{\mu,p})$ is isometric to $X$.
\end{theorem}

It is clear that taking into account this representation for any finite-dimensional subspace of $L_p$, the problem 
of embedding $\ell_2^k$ in subspaces of $L_p$ reduces to spaces $(\mathbb R^n, \|\cdot\|_{\mu, p})$ with $\mu$ satisfying
the condition \eqref{eq:iso}. Hence, the next paragraph is devoted to the study of these measures.

%%%%%%%%%%%%%%%%%%%%%%%%%%%%%%%%%%%%%%%%%%%%%%
\subsection{Isotropic measures on the sphere}
%%%%%%%%%%%%%%%%%%%%%%%%%%%%%%%%%%%%%%%%%%%%%%

An even Borel measure $\mu$ on $S^{n-1}$ is said to be {\it isotropic} if it
satisfies the following condition:
\begin{align*}
\|x\|_2^2 = \int_{S^{n-1}} |\langle x, \theta \rangle|^2 \, d\mu(\theta) ,
\end{align*} for all $x\in \mathbb R^n$. Equivalently, for all linear transformations $T:\mathbb R^n\to \mathbb R^n$ we have
\begin{align*}
{\rm trace}(T) =\int_{S^{n-1}} \langle \theta, T\theta \rangle \, d\mu(\theta).
\end{align*} For any such measure we may define 
the following family of centrally symmetric convex bodies $B_p(\mu)$ associated with $\mu$ and corresponding norms: 
\begin{align*}
x\mapsto \|x\|_{B_p(\mu)}:=\|\langle x, \cdot \rangle\|_{L_p(\mu)} = \left( \int_{S^{n-1}} |\langle x, z\rangle|^p \, d\mu(z)\right)^{1/p}, 
\quad 1\leq p <\infty .
\end{align*} The corresponding spaces, whose unit ball is $B_p(\mu)$, will be denoted by $X_p(\mu)$. Under this terminology and
notation, Lewis' theorem reads as follows:

\begin{theorem}[Lewis] \label{thm:Lewis}
Let $1\leq p<\infty$ and let $X$ be an $n$-dimensional subspace of $L_p$. Then, there exists an isotropic Borel 
measure $\mu$ on $S^{n-1}$ and a linear isometry $U:X_p(\mu) \to X$.
\end{theorem}

The next simple lemma collects several properties for the bodies $B_p(\mu)$.

\begin{lemma}\label{lem:meas-prop}
Let $\mu$ be a Borel isotropic measure on $S^{n-1}$ and let $Z$ be an $n$-dimensional standard Gaussian vector. Then, we 
have the following properties:
\begin{itemize}
\item [\rm i.] $\mathbb E \|Z\|_{B_q(\mu)}^q = \sigma_q^q \mu(S^{n-1})$, for $0<q<\infty$.

\item [\rm ii.] $\mu(S^{n-1})=n$.

\item [\rm iii.] For $p\geq 2$ we have $\|x\|_{B_p(\mu)} \leq \|x\|_2$ and for $1\leq p<q<\infty$ we have
$\|x\|_{B_p(\mu)} \leq n^{1/p-1/q} \|x\|_{B_q(\mu)}$, for all $x\in \mathbb R^n$.

\item [\rm iv.] (K. Ball) For every $1\leq p<\infty$ we have $|B_p(\mu)|\leq |B_p^n|$.

\item [\rm v.] For the body $B_q(\mu), \; q\geq 1$ we have $k(B_q(\mu))\geq cn^{\min\{1, 2/q\} }$.

\item [\rm vi.] There exists an absolute constant $c>0$ such that for all $2\leq q\leq c\log n$, 
one has $(\mathbb E \|Z\|_{B_q(\mu)}^2)^{1/2} \simeq q^{1/2}n^{1/q}$. 
In particular, for those $q$'s one has $k(B_q(\mu))\geq c qn^{2/q}$.
\end{itemize}
\end{lemma}

\noindent {\it Proof.} (i) We use Fubini's theorem and the rotation invariance of the Gaussian measure
to write
\begin{align*}
\mathbb E \|Z\|_{B_q(\mu)}^q=\int_{\mathbb R^n} \|x\|_{B_q(\mu)}^q \, d\gamma_n(x) = \int_{S^{n-1}} \int_{\mathbb R^n} |\langle x,\theta\rangle|^q \, d\gamma_n(x) \, d\mu(\theta) = \sigma_q^q \mu(S^{n-1}).
\end{align*} 
(ii) It follows from the above formula applied for $q=2$ and by employing the isotropic condition. 

\smallskip

\noindent (iii) Let $p\geq 2$. Note that
for all $u\in S^{n-1}$ we have
\begin{align*}
\|u\|_{B_p(\mu)}^p = \int_{S^{n-1}} |\langle u,\theta \rangle|^p \, d\mu(\theta) \leq \int_{S^{n-1}} |\langle u,\theta \rangle|^2\, d\mu(\theta) =1 .
\end{align*} For $1\leq p\leq q$ we apply H\"{o}lder's inequality
\begin{align*}
\|x\|_{B_p(\mu)} =\left( \int_{S^{n-1}} |\langle x, \theta\rangle|^p \, d\mu(\theta) \right)^{1/p} \leq \mu(S^{n-1})^{\frac{1}{p}-\frac{1}{q} } \left( \int_{S^{n-1}} |\langle x, \theta \rangle|^q \, d\mu(\theta) \right)^{1/q}.
\end{align*} 

\smallskip

\noindent (iv) This result was essentially proved by K. Ball in \cite{Ba}. A sketch of his very elegant proof is reproduced below 
for the sake of completeness.
Without loss of generality we may assume that $\mu$ is discrete, i.e. $\mu=\sum_{i=1}^m c_i \delta_{u_i}$, 
for some vectors $(u_i)$ in $S^{n-1}$ and positive numbers $(c_i)$ with $I=\sum_{i=1}^m  c_i u_i\otimes u_i$. Now we use the formula, 
which holds true for any centrally symmetric convex body $K$ on $\mathbb R^n$,
\begin{align*}
|K| = (\Gamma(1+n/p))^{-1} \int_{\mathbb R^n} e^{-\|z\|_K^p}\, dz,
\end{align*}  to get
\begin{align*}
|B_p(\mu)| = \frac{1}{\Gamma (1+\frac{n}{p})} \int_{\mathbb R^n} \prod_{i=1}^m f_i(\langle z,u_i\rangle) ^{c_i} \, dz,
\end{align*} where $f_i(t)=\exp(-|t|^p)$. The result follows by the Brascamp-Lieb inequality. 

\smallskip

\noindent (v) First consider the case $2<q<\infty$. For the critical dimension of 
the space $X_q(\mu)=(\mathbb R^n, \|\cdot \|_{B_q(\mu)})$, note 
that $k(X_q(\mu)) = \mathbb E \|Z\|_{B_q(\mu)}^2 / b^2( B_q(\mu) ) \geq n^{2/q}$ by the third assertion. 

\noindent Now we turn in the range $1\leq q\leq 2$. Using H\"{o}lder's inequality we may write
\begin{align*}
\left(\mathbb E \|Z\|_{B_q(\mu)}^2 \right)^{1/2} \geq n^{1/2} \left( \frac{|B_2^n|}{|B_q(\mu)|} \right)^{1/n} 
\geq n^{1/2} \left( \frac{|B_2^n|}{|B_q^n|} \right)^{1/n} \simeq n^{1/q} ,
\end{align*} where in the last step we have used Ball's volumetric estimate (iv). The result follows once 
we recall that $b(B_q(\mu)) \leq n^{1/q-1/2}$ for $1\leq q\leq 2$.

\smallskip

\noindent (vi) We define the parameter
\begin{align*} 
q_0\equiv q_0(\mu): = \max \left\{ q\in [2,n] : k(B_p(\mu))\geq p, \; \forall p\in [2,q] \right\}.
\end{align*} By the continuity of the map
$p\mapsto k(B_p(\mu))$ and the fact that $k(B_q(\mu))\leq n$ for all $q\geq 2$, while $k(B_2(\mu))=n$ we get $q_0= k(B_{q_0}(\mu))$. Lemma
\ref{lem:conc-log-sob} shows that 
$\left( \mathbb E \|Z\|_{B_{q_0}(\mu)}^{q_0} \right)^{1/q_0} \leq c_1 \left( \mathbb E \|Z\|_{B_{q_0}(\mu)}^2 \right)^{1/2}$, so we may write
\begin{align*}
q_0 =k (B_{q_0}(\mu)) =\frac{\mathbb E \|Z\|_{B_{q_0}(\mu)}^2}{b^2(B_{q_0}(\mu))}\geq c_1^{-2} (\mathbb E \|Z\|_{B_{q_0}(\mu)}^{q_0})^{2/q_0} =
c_1^{-2} \sigma_{q_0}^2 n^{2/q_0} \Longrightarrow q_0 \geq c_2 \log n.
\end{align*} Therefore, by the definition of $q_0$ we have $k(B_q(\mu))\geq q$ for all $2\leq q\leq q_0$ and by Lemma \ref{lem:conc-log-sob}
again, we get
\begin{align*}
\sigma_q n^{1/q}= \left(\mathbb E \|Z\|_{B_q(\mu)}^q\right)^{1/q} \leq c_1 \left( \mathbb E \|g\|_{B_q(\mu)}^2 \right)^{1/2}.
\end{align*}
Moreover, we have
\begin{align*}
k(B_q(\mu)) = \frac{\mathbb E \|Z\|_{B_q(\mu)}^2}{b^2(B_q(\mu))}\geq c_1^{-2} (\mathbb E \|Z\|_{B_q(\mu)}^q)^{2/q} =
c_1^{-2} \sigma_q^2 n^{2/q} \geq c_3 q n^{2/q}. \end{align*} This can be interpreted as $k(B_q(\mu))\geq ck(\ell_q^n)$, provided
that $2\leq q\leq c\log n$ for some absolute constant $c>0$. For a proof of the fact that $k(\ell_q^n)\simeq qn^{2/q}$ when $2\leq q\leq c\log n$
the reader is referred to \cite{SZin}. \prend

\begin{lemma}\label{lem:moms-meas}
Let $\mu$ be a Borel isotropic measure on $S^{n-1}$. For $q\geq 2$ and for all $r\geq 1$ we have
\begin{align*}
I_{rq}(\gamma_n,B_q(\mu))/I_q(\gamma_n, B_q(\mu)) \leq \sqrt{1+\frac{q(r-1)}{\sigma_q^2 n^{2/q} } } \leq \sqrt{1+\frac{c(r-1)}{n^{2/q} } },
\end{align*} where $c>0$ is an absolute constant.
\end{lemma}

\noindent {\it Proof.} Note that Lemma \ref{lem:meas-prop} (iii) implies 
$\left| \|x\|_{B_q(\mu)} - \|y\|_{B_q(\mu)} \right|\leq \|x-y\|_2$ for all $x, y\in \mathbb R^n$. Hence, if we use 
Lemma \ref{lem:conc-log-sob} we obtain
\begin{align*}
\left( \frac{I_{rq}}{I_q} \right)^2 \leq 1+ \frac{q(r-1)}{I_q^2} = 1+\frac{q(r-1)}{\sigma_q^2 n^{2/q}} ,
\end{align*} where the last estimate follows from Lemma \ref{lem:meas-prop}. Finally, using the fact that $\sigma_q\simeq \sqrt{q}$
we conclude the second estimate. \prend

%%%%%%%%%%%%%%%%%%%%%%%%%%%%%%%%%%%%%%%%%%%%%%%%%%
\subsection{A Gaussian inequality}
%%%%%%%%%%%%%%%%%%%%%%%%%%%%%%%%%%%%%%%%%%%%%%%%%%

The next inequality is due to Pisier (for a proof see \cite{Pis}).

\begin{theorem}\label{thm:Mau-Pis}
Let $\phi:\mathbb R\to \mathbb R$ be a convex function and let $f:\mathbb R^n\to \mathbb R$ be $C^1$-smooth. Then, if $Z,W$ are
independent copies of a Gaussian random vector, we have
\begin{align*} 
\mathbb E \phi\left( f(Z)-f(W)\right)\leq  \mathbb E \phi \left( \frac{\pi}{2} \langle \nabla f(Z), W\rangle \right).
\end{align*} 
\end{theorem}

Here we prove a generalization of this inequality in the context 
of Gaussian processes generated by the action of a random matrix with i.i.d standard Gaussian entries on a fixed vector in $S^{n-1}$.
The next inequality was stated in \cite{PVZ} without a proof. Below we give the details for reader's convenience.

\begin{theorem}\label{thm:Mau-Pis-stoch}
Let $\phi:\mathbb R\to \mathbb R$ be a convex function and let $f:\mathbb R^n\to \mathbb R$ be $C^1$-smooth. If 
$G=(g_{ij})_{i,j=1}^{n,k}$ is a Gaussian matrix and $a,b\in S^{k-1}$, then we have
\begin{align*} 
\mathbb E \phi\left( f(Ga)-f(Gb)\right)\leq  \mathbb E \phi \left( \frac{\pi}{2} \|a-b\|_2 \langle \nabla f(Z), W\rangle \right),
\end{align*} where $Z,W$ are independent copies of a standard Gaussian $n$-dimensional random vector. 
\end{theorem} 

\noindent {\it Proof.} If $a=b$ then, there is nothing to prove. If $a=-b$ then, by setting $F(z)=f(z)-f(-z)$ we may write:
\begin{align*}
\mathbb E \phi\left( f(Ga)-f(Gb)\right) = \mathbb E \phi( F(Z) ) \leq \mathbb E \phi (F(Z)-F(W)),
\end{align*} for $Z,W$ independent copies of a standard Gaussian random vector, where we have used the fact $\mathbb EF(Z)=0$ and
Jensen's inequality. Then, a direct application of Theorem \ref{thm:Mau-Pis} yields:
\begin{align*}
\mathbb E \phi (F(Z)-F(W)) &\leq \mathbb E \phi\left( \frac{\pi \langle \nabla f(Z), W \rangle + \pi \langle \nabla f(-Z),W \rangle}{2}\right) \\
& \leq \mathbb E \frac{ \phi(\pi \langle \nabla f(Z), W \rangle) + \phi(\pi \langle \nabla f(-Z),W \rangle)}{2} \\
&= \mathbb E \phi( \pi \langle \nabla f(Z),W\rangle),
\end{align*} by the convexity of $\phi$.

In the general case, fix $a,b\in S^{k-1}$ with $a\neq \pm b$ and define $p:=\frac{a+b}{2}$. Note that since $\|a\|_2=\|b\|_2$ we have that the
vector $u:=a-p$ is perpendicular to $p$. Set $W:=G(u)$ and $Z:=G(p)$ and note that $W, Z$ are independent random vectors in $\mathbb R^n$ with
$W \sim N({\bf 0}, \|u\|_2^2I_n)$, $Z \sim N({\bf 0}, \|p\|_2^2I_n)$. Since $G(a)=Z+W$ and $G(b)=Z-W$, we may write:
\begin{align*}
\mathbb E \phi\left( f(Ga) -f(Gb) \right) =\mathbb E_Z\mathbb E_W \phi \left( f(Z+W) -f(Z-W) \right).
\end{align*}
Denote $F(w,z):=f(z+w)-f(z-w)$. Then, we may write:
\begin{align*} 
\mathbb E \phi(f(Ga) -f(Gb)) =\iint \phi( F(w,z)) \, d\gamma_{n,\sigma_1}(w) \, d\gamma_{n,\sigma_2}(z) ,
\end{align*} where $\sigma_1=\|u\|_2>0, \, \sigma_2=\|p\|_2>0$. For fixed $z$, we may apply Theorem \ref{thm:Mau-Pis} to the function 
$w\mapsto F(w,z)$ (note that $\int F(w,z) \, d\gamma_{n, \sigma_1}(w)=0$) to get:
\begin{align*}
\int \phi( F(w,z) )\, d\gamma_{n,\sigma_1}(w) &\leq \iint \phi\left(\frac{\pi}{2} \langle \nabla_w F(w,z),y \rangle \right)\, d\gamma_{n,\sigma_1}(w) \, d\gamma_{n,\sigma_1}(y) \\
&\leq \iint \frac{\phi(\pi \langle \nabla f(w+z),y\rangle ) + \phi(\pi \langle \nabla f(z-w),y\rangle) }{2} \, d\gamma_{n,\sigma_1}(w) \, d\gamma_{n,\sigma_1}(y) \\
& = \iint \phi\left( \pi \langle \nabla f(w+z),y \rangle \right) \,d\gamma_{n,\sigma_1}(w) \, d\gamma_{n,\sigma_1}(y),
\end{align*} by the convexity of $\phi$. Integration with respect to $\gamma_{n,\sigma_2}$ over $z$ provides:
\begin{align*}
\iint \phi( F(w,z)) \, d\gamma_{n,\sigma_1}(w) d\gamma_{n,\sigma_2}(z) 
&\leq \int \left[ \iint \phi \left( \pi \langle \nabla f(w+z),y \rangle \right) d\gamma_{n,\sigma_1}(w) d\gamma_{n,\sigma_2}(z)\right] d\gamma_{n,\sigma_1}(y) \\
&= \int \left[ \int \phi\left( \pi \langle \nabla f(x),y \rangle\right) d(\gamma_{n,\sigma_1} \ast \gamma_{n,\sigma_2})(x)\right] d\gamma_{n,\sigma_1}(y) \\
& = \iint \phi \left( \pi \sigma_1 \langle \nabla f(x),y \rangle \right) d\gamma_n(x) \, d\gamma_n(y), 
\end{align*} where we have used the fact that $\gamma_{n,\sigma_1} \ast \gamma_{n,\sigma_2} =\gamma_{n,\sigma_1^2+\sigma_2^2} \equiv \gamma_n$, 
since $\sigma_1^2+\sigma_2^2=\|a\|_2^2=1$. The result follows. \prend

\begin{remark}  \label{rem:2.8.1}\rm 1. Applying this for $\phi(t)=|t|^r, \; r\geq 1$ and taking into account the invariance of the Gaussian measure under orthogonal 
transformations we derive the next $(r,r)$-Poincar\'{e} inequalities:
\begin{align}\label{eq: 2.23}
\left( \mathbb E | f(Ga)-f(Gb)|^r \right)^{1/r} \leq C \sqrt{r} \|a-b\|_2 \left(\mathbb E \|\nabla f(Z)\|_2^r \right)^{1/r},
\end{align} for $a,b\in S^{k-1}$, where $Z$ is a standard Gaussian random vector in $\mathbb R^n$.

\noindent 2. Assuming further that $f$ is $L$-Lipschitz we may apply Theorem \ref{thm:Mau-Pis-stoch} 
for $\phi(t)=e^{\lambda t}, \; \lambda>0$ to get:
\begin{align} \label{eq: 2.27}
\mathbb E \exp \left(\lambda \left(f(Ga)-f(Gb) \right) \right) \leq \mathbb E \exp \left(\lambda^2 \frac{\pi^2}{2}\|a-b\|_2^2 \|\nabla f(Z)\|_2^2 \right) \leq 
\exp \left(\lambda^2 \frac{\pi^2}{2}\|a-b\|_2^2 L^2\right).
\end{align}
Then Markov's inequality yields Schechtman's distributional
inequality from \cite{Sch1}:
\begin{align} \label{eq: 2.24}
\mathbb P \left( |f(Ga)-f(Gb)| > t \right)\leq C\exp \left( -c t^2 / ( \|a-b\|^2 L^2)\right),
\end{align} for all $t>0$, where $a,b\in S^{k-1}$. Let us note that \eqref{eq: 2.27} for $f$ being a norm, has also appeared in \cite{Schmu}.

\smallskip

\noindent 3. For $a,b\in S^{k-1}$ with $\langle a,b\rangle=0$ the matrix $G$ generates the vectors $Z=Ga$ and $W=Gb$ which
are independent copies of a standard $n$-dimensional Gaussian random vector. For example, inequality \eqref{eq: 2.24} reduces
to the classical concentration inequality:
\begin{align} \label{eq: Gauss-conc}
\mathbb P \left( | f(Z) - f(W) | > t \right) \leq C\exp \left( -c t^2 /  L^2 \right),
\end{align} for all $t>0$.
\end{remark}

%%%%%%%%%%%%%%%%%%%%%%%%%%%%%%%%%%%%%%%%%%
\section{Gaussian concentration for $B_p(\mu)$ norms}
%%%%%%%%%%%%%%%%%%%%%%%%%%%%%%%%%%%%%%%%%%

A direct application of the Gaussian concentration inequality \eqref{eq: Gauss-conc} for the norms $\|\cdot\|_{B_p(\mu)}$, $2<p<\infty$ implies:
\begin{align} \label{eq:3.1}
\mathbb P \left( \big| \|Z\|_{B_p(\mu)} -I_1\big| >t I_1\right) \leq C \exp(-ct^2 I_1^2) \leq C\exp(-ct^2 n^{2/p}),
\end{align} for all $t>0$, where $I_1\equiv I_1(\gamma_n,B_p(\mu))$. It is known (see \cite{PVZ}) that the large
deviation estimate ($t\geq 1$) the inequality \eqref{eq:3.1} provides 
is sharp (up to constants).

In this paragraph we prove that for $2< p <\infty$ and $\mu$ isotropic Borel measure on $S^{n-1}$, the bodies $B_p(\mu)$ exhibit better
concentration ($0<t<1$) than the one implied by the Gaussian concentration inequality on $\mathbb R^n$ in terms of the Lipschitz constant. 
Later, this will be used to prove the announced dependence on $\varepsilon$ and $n$ in Dvoretzky's theorem for any $n$-dimensional 
subspace of $L_p$. Our main tool is the probabilistic inequality proved in Theorem \ref{thm:Mau-Pis-stoch} and as was formulated further 
in Remark \ref{rem:2.8.1}.1.

We apply inequality \eqref{eq: 2.23} for $f(x)=\|x\|^p= \int |\langle x,\theta \rangle|^p \, d\mu(\theta)$. To this end
we have to compute the gradient. Note that
\begin{align*}
\|\nabla f (x) \|_2^2 = p^2 \sum_{i=1}^n \left| \int_{S^{n-1}} \theta_i | \langle x,\theta \rangle|^{p-1} \sgn{\langle x,\theta \rangle} \, d\mu(\theta) \right|^2.
\end{align*} We also have the following:

\smallskip

\noindent {\it Claim.} For almost every $x\in\mathbb R^n$ we have
\begin{align*}
\|\nabla f(x)\|_2^2 \leq p^2 \|x\|_{B_{2p-2}(\mu)}^{2p-2}.
\end{align*} 

\noindent {\it Proof of Claim.} Let $b_i\equiv b_i(x): =\int_{S^{n-1}} |\langle x, z\rangle |^{p-1} \sgn{\langle x, z \rangle} z_i \, d\mu(z)$. 
Using duality we may write
\begin{align*}
\sum_{i=1}^n \left(\int_{S^{n-1}} |\langle x, z\rangle |^{p-1} \sgn{\langle x, z \rangle} z_i \, d\mu(z) \right)^2 &=
\max_{\theta \in S^{n-1}} \left| \sum_{i=1}^n b_i \theta_i\right|^2 \\
&= \max_{\theta \in S^{n-1}} \left| \int_{S^{n-1}} |\langle x, z\rangle|^{p-1} \sgn{\langle x, z\rangle} \langle z,\theta \rangle \, d\mu(z) \right|^2 \\
&\leq \int_{S^{n-1}} |\langle x, z\rangle|^{2p-2}\, d\mu(z),
\end{align*} where we have used the Cauchy-Schwarz inequality and the isotropic condition. \prend

\smallskip

\noindent Therefore, using the Claim and the inequality \eqref{eq: 2.23} we get for every $a,b\in S^{k-1}$
\begin{align*}
\left( \mathbb E |f(Ga)-f(Gb)|^r \right)^{1/r} \leq C p r^{1/2} \|a-b\|_2  \left( \mathbb E \|Z\|_{B_{2p-2}(\mu)}^{r(p-1)} \right)^{1/r},
\end{align*} for all $r\geq 1$. By employing Lemma \ref{lem:moms-meas} we find
\begin{align*}
\left( \mathbb E |f(Ga)-f(Gb)|^r \right)^{1/r} 
&\leq 
C p r^{1/2} \|a-b\|_2 \left( \mathbb E\|Z\|_{B_{2p-2}(\mu)}^{2p-2}\right)^{1/2} \left(1+\frac{(r-2)(p-1)}{\sigma_{2p-2}^2 n^{\frac{1}{p-1}}} \right)^{\frac{p-1}{2}}\\
& < C p r^{1/2} \|a-b\|_2 \sigma_{2p-2}^{p-1} n^{1/2} \left(1+\frac{r(p-1)}{\sigma_{2p-2}^2 n^{\frac{1}{p-1}}} \right)^{\frac{p-1}{2}} \\
& < C p \|a-b\|_2 \sigma_{2p-2}^{p-1} n^{1/2} 2^{\frac{p-1}{2}} \max\left\{ r^{1/2}, \frac{r^{p/2} (p-1)^{\frac{p-1}{2}} }{\sigma_{2p-2}^{p-1} n^{1/2}}\right\},
\end{align*} for all $r\geq 2$. We define
\begin{align} \label{eq: def-a}
\alpha(n,p,r): =\max\left\{ r^{1/2}, \frac{r^{p/2} (p-1)^{\frac{p-1}{2}} }{\sigma_{2p-2}^{p-1} n^{1/2}}\right\} , \quad r > 0
\end{align} and we summarize the above discussion to the following:

\begin{proposition} \label{prop:main-1}
Let $2< p < \infty $ and let $\mu$ be a Borel isotropic measure on $S^{n-1}$. If $G=(g_{ij})_{i,j=1}^{n,k}$ is a
Gaussian matrix and $a,b\in S^{k-1}$, then we have
\begin{align*}
\left( \mathbb E \big| \|Ga\|_{B_p(\mu)}^p - \|Gb \|_{B_p(\mu)}^p \big|^r \right)^{1/r} \leq 
C p \|a-b\|_2\sigma_{2p-2}^{p-1} n^{1/2} 2^{\frac{p}{2}} \alpha(n,p,r),
\end{align*} for all $r\geq 2$, where $\alpha(n,p,\cdot)$ is defined in \eqref{eq: def-a}
\end{proposition}

We are now ready to prove the main result of this section.

\begin{theorem} \label{thm:conc-B_p}
Let $2<p<\infty$ and let $\mu$ be a Borel isotropic measure on $S^{n-1}$ with $n>e^p$. Then, we have
\begin{align*}
\mathbb P \left( \big| \|Z\|_{B_p(\mu)} - ( \mathbb E\|Z\|_{B_p(\mu)}^p)^{1/p} \big| \geq \varepsilon (\mathbb E \|Z\|_{B_p(\mu)}^p)^{1/p} \right) \leq C \exp \left( -c \psi(n,p,\varepsilon) \right),
\end{align*} for every $\varepsilon>0$, where $\psi(n,p, \cdot)$ is defined by
\begin{align} \label{eq: def-psi}
\psi(n,p,t) := \min\left\{ \frac{t^2 n}{p4^p}, (tn)^{2/p} \right\}, \, t>0,
\end{align} and $C,c>0$ are absolute constants.
\end{theorem}

\noindent {\it Proof.} Using Proposition \ref{prop:main-1} for $a,b\in S^{k-1}$ with $\langle a,b\rangle=0$ and applying 
Jensen's inequality we obtain
\begin{align*}
\left( \mathbb E \big| \|Z\|_{B_p(\mu)}^p - \mathbb E \|Z \|_{B_p(\mu)}^p \big|^r \right)^{1/r} \leq 
C p \sigma_{2p-2}^{p-1} n^{1/2} 2^{p/2} \alpha(n,p,r),
\end{align*} for all $r\geq 2$. Therefore Markov's inequality yields
\begin{align*}
\mathbb P \left( \big| \|Z\|_{B_p(\mu)}^p - \mathbb E \|Z\|_{B_p(\mu)}^p \big| >\varepsilon \right) \leq \left( \frac{C p \sigma_{2p-2}^{p-1} n^{1/2} 2^{p/2} \alpha(n,p,r) }{\varepsilon} \right)^r.
\end{align*} Note that the inverse of the map $r\mapsto \alpha(n,p,r)$ is given by 
\begin{align*}
\alpha^{-1}(n,p,s) = \min \left\{ s^2, \, \frac{ s^{2/p} n^{1/p} \sigma_{2p-2}^{\frac{2p-2}{p} }  }{ (p-1)^{\frac{p-1}{p} } } \right\}, \quad s>0,
\end{align*} thus we may choose $r_\varepsilon \geq 2$ such that $\alpha(n,p,r_\varepsilon) =\frac{\varepsilon}{ e C p \sigma_{2p-2}^{p-1} n^{1/2} 2^{p/2}}$, as long as 
the range of $\varepsilon>0$ satisfies $ \alpha(n,p,r_\varepsilon) \geq \alpha(n,p,2) $. Otherwise $\alpha(n,p,r_\varepsilon)< \alpha(n,p,2)\simeq \max\{1, (e^p/n)^{1/2}\}
\simeq 1$ provided that $n$ is large enough with respect to $p$. Hence, we get
\begin{align*}
\mathbb P \left( \big| \|Z\|_{B_p(\mu)}^p - \mathbb E \|Z\|_{B_p(\mu)}^p \big| >\varepsilon \right) \leq C_1 \exp \left( -\alpha^{-1} \left( n,p, \frac{\varepsilon}{ e C p \sigma_{2p-2}^{p-1} n^{1/2} 2^{p/2}} \right) \right),
\end{align*} for all $\varepsilon>0$. We may check that
\begin{align*}
\alpha^{-1} \left( n,p, \frac{\varepsilon}{ e C p \sigma_{2p-2}^{p-1} n^{1/2} 2^{p/2} } \right) \simeq 
\min \left\{ \frac{\varepsilon^2}{n p^2 2^p \sigma_{2p-2}^{2p-2} } , \; \frac{\varepsilon^{2/p} }{p} \right\},
\end{align*} which implies
\begin{align*}
\mathbb P \left( \big| \|Z\|_{B_p(\mu)}^p - \mathbb E \|Z\|_{B_p(\mu)}^p \big| >\varepsilon \right) \leq 
C_1 \exp \left( -c_1 \min \left\{ \frac{\varepsilon^2}{ n p^2 2^p \sigma_{2p-2}^{2p-2}} , \; \frac{\varepsilon^{2/p} }{p} \right\} \right),
\end{align*} for every $\varepsilon>0 $. It follows that
\begin{align*}
\mathbb P \left( \big| \|Z\|_{B_p(\mu)}^p - \mathbb E \|Z\|_{B_p(\mu)}^p \big| >\varepsilon \mathbb E \|Z\|_{B_p(\mu)}^p \right) \leq 
C_1 \exp \left( -c_1 \min \left\{ \frac{\varepsilon^2 n \sigma_p^{2p}}{ p^2 2^p \sigma_{2p-2}^{2p-2}} , \; \frac{(\varepsilon n)^{2/p} \sigma_p^{2} }{p} \right\} \right),
\end{align*} for every $\varepsilon>0$. The asymptotic estimate \eqref{eq: 2.3} yields $\sigma_p^{2p} / \sigma_{2p-2}^{2p-2} \simeq p2^{-p}$ 
and $\sigma_p \simeq p^{1/2}$, thus we conclude
\begin{align} \label{eq: 3.15}
\mathbb P \left( \big| \|Z\|_{B_p(\mu)}^p - \mathbb E \|Z\|_{B_p(\mu)}^p \big| >\varepsilon \mathbb E \|Z\|_{B_p(\mu)}^p \right) \leq 
C_1 \exp \left( -c_1' \min \left\{ \frac{\varepsilon^2 n }{ p 4^p } , \; (\varepsilon n)^{2/p}  \right\} \right),
\end{align} for all $\varepsilon>0$. This further implies that
\begin{align*}
\mathbb P \left( \big| \|Z\|_{B_p(\mu)} - \left( \mathbb E \|Z\|_{B_p(\mu)}^p \right)^{1/p} \big| >\varepsilon \left( \mathbb E \|Z\|_{B_p(\mu)}^p \right)^{1/p} \right) 
\leq 2C_1 \exp \left( -c_1' \min \left\{ \frac{\varepsilon^2 n }{ p 4^p } , \; (\varepsilon n)^{2/p}  \right\} \right),
\end{align*} for all $\varepsilon>0$.  In order to verify the latter we may write
\begin{align*}
\mathbb P \left( \|Z\|_{B_p(\mu)}  > (1+\varepsilon) \left( \mathbb E \|Z\|_{B_p(\mu)}^p \right)^{1/p} \right) &\leq 
\mathbb P \left( \|Z\|_{B_p(\mu)}^p  > (1+\varepsilon) \mathbb E \|Z\|_{B_p(\mu)}^p  \right) \\
&\leq 
C_1 \exp \left( -c_1' \min \left\{ \frac{\varepsilon^2 n }{ p 4^p } , \; (\varepsilon n)^{2/p}  \right\} \right) ,
\end{align*} for all $\varepsilon>0$ by the estimate \eqref{eq: 3.15}. We argue similarly for the other case. \prend

\medskip

\begin{remark} By the well known symmetrization argument for any random variable $\xi$
\begin{align*}
\mathbb P ( |\xi -{\rm med}(\xi)|>t ) \leq 4 \inf_{\alpha\in \mathbb R} \mathbb P(|\xi-\alpha| >t/2), \quad  t>0, 
\end{align*} we may replace $(\mathbb E\|Z\|_{B_p(\mu)}^p)^{1/p}$ by a 
median of $x\mapsto \|x\|_{B_p(\mu)}$ (or the expected value $\mathbb E \|X\|_{B_p(\mu)}$) with respect to the Gaussian
measure $\gamma_n$ (see also \cite[Appendix V]{MS}).
\end{remark}

\medskip

We shall also need the next variant of Theorem \ref{thm:conc-B_p}.

\begin{theorem} \label{thm:conc-B_p-2}
Let $2<p<\infty$ and let $\mu$ be a Borel isotropic probability measure on $S^{n-1}$ with $n>e^p$. If $G=(g_{ij})_{i,j=1}^{n,k}$ is a 
Gaussian matrix and $a,b \in S^{k-1}$, then
\begin{align*}
\mathbb P \left( \big| \|Ga\|_{B_p(\mu)}^p - \|Gb\|_{B_p(\mu)}^p \big| > t \mathbb E\|Z\|_{B_p(\mu)}^p \right) \leq
 C\exp \left(-c \psi \left(n,p, \frac{t}{\|a-b\|_2} \right) \right),
\end{align*} for all $t>0$, where $\psi(n,p,\cdot)$ is defined in \eqref{eq: def-psi}.
\end{theorem}

\noindent {\it Proof.} The proof is similar to the proof of Theorem \ref{thm:conc-B_p}. We omit the details. \prend

\medskip

The following estimate on the fluctuations of the norm $x\mapsto \|x\|_{B_p(\mu)}$ is immediate:

\begin{corollary} [Gaussian variance for $B_p(\mu)$] \label{cor:var-B_p}
Let $\mu$ be an isotropic Borel measure on $S^{n-1}$ and let $1\leq p<\infty$. Then, 
\begin{align*}
{\rm Var} \|Z\|_{B_p(\mu)} \leq e^{cp} n^{2/p-1}.
\end{align*} In particular, we have
\begin{align*}
\frac{{\rm Var}\|Z\|_{B_p(\mu)}}{\mathbb E\|Z\|^2_{B_p(\mu)}} \leq \frac{e^{cp} }{n},
\end{align*} where $Z$ is a standard Gaussian vector.
\end{corollary}

\noindent {\it Proof.} We distinguish two cases. For $1\leq p \leq 2$ we may bound as follows:
\begin{align*}
{\rm Var}\|Z\|_{B_p(\mu)} \lesssim b^2(B_p(\mu)) \lesssim n^{2/p-1},
\end{align*} where we have used Lemma  \ref{lem:meas-prop}. For $2<p<\infty$ consider $Z'$ an independet copy of $Z$ to write
\begin{align*}
2 {\rm Var} \|Z\|_{B_p(\mu)} =\mathbb E (\|Z\|_{B_p(\mu)} -\|Z'\|_{B_p(\mu)})^2 \leq 
\frac{1}{p^2} \mathbb E \left( \frac{\|Z\|_{B_p(\mu)}^p -\|Z'\|_{B_p(\mu)}^p }{ \min\{ \|Z\|_{B_p(\mu)}^{p-1}, \| Z' \|_{B_p(\mu)}^{p-1}  \} } \right)^2,
\end{align*} where we have used the numerical inequality $|a^p -b^p| \geq p |a-b| \min\{ a^{p-1}, b^{p-1}\}$ for $a,b>0$ and $p>1$. 
The Cauchy-Schwarz inequality implies
\begin{align*}
{ \rm Var} \|Z\|_{B_p(\mu)} \leq \frac{1}{p^2} \frac{ \left( \mathbb E \left| \|Z\|_{B_p(\mu)}^p - \|Z' \|_{B_p(\mu)}^p \right|^4 \right)^{1/2} }{ I_{-4(p-1)}^{2(p-1)} (\gamma_n, B_p(\mu))}.
\end{align*} The numerator is directly estimated by Proposition \ref{prop:main-1}. For the denominator we employ the main result of \cite{KV} along
with the fact $k(B_p(\mu)) \geq c_1pn^{2/p}$ for $n\geq e^{C_1p}$ (proved in Lemma  \ref{lem:meas-prop}.vi) to obtain
\[ I_{-4(p-1)}(\gamma_n, B_p(\mu)) \geq c_2  I_p(\gamma_n, B_p(\mu)) = \sigma_p n^{1/p},
\] by Lemma \ref{lem:meas-prop}.i. Combining all the above we arrive at the desired estimate. The details are left to the reader. \prend

%%%%%%%%%%%%%%%%%%%%%%%%%%%%%%%%%%%%%%%%%%%%%%%%%%%%%%%%%%
\section{Embedding $\ell_2^k$ into subspaces of $L_p$ for $2<p<\infty$}
%%%%%%%%%%%%%%%%%%%%%%%%%%%%%%%%%%%%%%%%%%%%%%%%%%%%%%%%%%

In this paragraph we prove the improved estimate on Dvoretzky's theorem for the subspaces of $L_p, \; 2<p<\infty$.

\begin{theorem} \label{thm:main-2}
Let $2<p<\infty$. Then for every $n$-dimensional subspace $X$ of $L_p$ and any $0<\varepsilon<1$ there exists 
$k\geq c_p \psi(n,p,\varepsilon)$ and linear map
$T: \ell_2^k\to X$ such that $\|x\|_2 \leq \|Tx\|_X \leq (1+\varepsilon) \|x\|_2$ for
all $x\in \ell_2^k$, where $c_p>0$ is constant depending only on $p$ and $\psi(n,p,\cdot)$ is given by \eqref{eq: def-psi}.
\end{theorem}

Clearly, the above theorem follows from the random version of Dvoretzky's theorem for subspaces of $L_p$ whose unit ball is
in Lewis' position, or equivalently for the bodies $B_p(\mu)$ with $\mu$ isotropic measure on $S^{n-1}$. More precisely, we have the following:

\begin{theorem} \label{thm:main-3}
Let $2<p<\infty$ and let $X$ be an $n$-dimensional subspace of $L_p$ represented on $\mathbb R^n$ 
whose unit ball $B_X$ is in Lewis' position. Then, for any
$\varepsilon \in (0,1)$ there exists $k\geq c_p\psi(n,p,\varepsilon)$ such that the random $k$-dimensional subspace $F$ of $X$ is 
$(1+\varepsilon)$-spherical with probability greater than $1-e^{-c_p\psi(n,p,\varepsilon)}$, 
where $c_p>0$ depends only on $p$ and $\psi(n,p,\cdot)$ is defined in \eqref{eq: def-psi}.
\end{theorem}

Let us note that once we have established 
the concentration estimate of Theorem \ref{thm:conc-B_p}, then a standard net argument yields the result with an extra 
$\log(1/\varepsilon)$ term. Indeed; fix $k\leq n$ and let $G=(g_{ij})_{i,j=1}^{n,k}$ be a Gaussian matrix with independent standard
entries. Let ${\cal D}$ be a $\delta$-net on $S^{k-1}$ with cardinality $| {\cal D}| \leq (3/\delta)^k$ (see \cite[Lemma 2.6]{MS} for the details). Then using the union bound, Theorem \ref{thm:conc-B_p} and the fact that $Gu$ is equidistributed to $Z \sim N({\bf 0},I_n)$ 
we obtain
\begin{align*}
\mathbb P \left( \exists \, u\in {\cal D} : \left| \|Gu\|_{B_p(\mu)} - \mathbb E\|Z\|_{B_p(\mu)} \right|  \geq \varepsilon \mathbb E\|Z\|_{B_p(\mu)} \right) \leq 
(3 / \delta)^k C \exp \left( -c\psi(n,p,\varepsilon) \right).
\end{align*} Choosing $\delta\simeq \varepsilon$ we find that with probability greater than $1-e^{-c'\psi(n,p,\varepsilon)}$ the random operator
$G$ satisfies
\[ \left| \|Gu\|_{B_p(\mu)} - \mathbb E\|Z\|_{B_p(\mu)} \right|  \leq \varepsilon \mathbb E\|Z\|_{B_p(\mu)} , \] 
for all $u\in {\cal D}$, as long as $k\lesssim (\log\frac{1}{\varepsilon})^{-1} \psi(n,p,\varepsilon)$. It's routine to check that 
we may pass to the whole sphere $S^{k-1}$ at cost of an oscillation at most $2\varepsilon \mathbb E\|Z\|_{B_p(\mu)}$, 
see e.g. \cite[Lemma 4.1]{MS}. 

Theorem \ref{thm:conc-B_p-2} serves exactly the purpose of removing this term. Then we use this inequality along with a chaining method 
to conclude the logarithmic-free dependence on $\varepsilon$ in our main result. This approach has been inspired by \cite{Sch1}. 
However, the method from \cite{Sch1} is not directly applicable here, since it lies in estimates involving the Lipschitz constant. 
As we have already explained
such estimates would only yield suboptimal bounds and one has to keep track of the higher moments of the length of the gradient until the very 
last step. This forces us to establish the inequality in Theorem \ref{thm:Mau-Pis-stoch}. In probabilistic terms Theorem \ref{thm:conc-B_p-2} says that the process $(\|G\theta\|_{B_p(\mu)}^p-I_p^p)_{\theta \in S^{k-1} }$ has two-level tail behavior described by $\psi(n,p,\cdot)$.

\medskip

Now we turn to proving the main result.

\medskip

\noindent {\it Proof of Theorem  \ref{thm:main-3}. } Let $2<p<\infty$ and let $X$ be an $n$-dimensional subspace of $L_p$ whose 
unit ball is in Lewis' position. Then, 
Lewis' theorem (Theorem \ref{thm:Lewis}) yields the existence of an isotropic Borel measure $\mu$ on $S^{n-1}$ and a linear 
isometry $S: X_p(\mu) \to X$, hence we may identify $X$ with $X_p(\mu)$.
We have to show that the ball $B_p(\mu)$ has random almost spherical $k$-dimensional sections 
with $k$ as large as possible. Let $\{g_{ij}(\omega) \}_{i,j=1}^{n,k}$ be i.i.d. standard
normals in some probability space $(\Omega, P)$ and consider the random Gaussian operator 
$G_\omega=(g_{ij}(\omega))_{i,j=1}^{n,k} :\ell_2^k \to X_p(\mu)$.
We will prove that with overwhelming probability the operator $G$ is $(1+\varepsilon)$-isomorphic embedding when $k$ is sufficiently large.
Toward this end, we employ Theorem \ref{thm:conc-B_p-2} and the chaining argument from \cite{Sch1}.
For each $j=1,2,\ldots$ consider $\delta_j$-nets ${\cal N}_j$ on $S^{k-1}$ with cardinality $|{\cal N}_j| \leq (3/\delta_j)^k$ 
(see \cite[Lemma 2.6]{MS}). Note that for any
$\theta\in S^{k-1}$ and for  all $j$ there exist $u_j\in {\cal N}_j$ with $\|\theta-u_j\|_2\leq \delta_j$ and by the triangle inequality it follows that
$\|u_j-u_{j-1}\|_2\leq \delta_j+\delta_{j-1}$. Moreover, if we assume that $\delta_j\to 0$ as $j\to \infty$ and $(t_j)$ is a sequence of numbers 
with $t_j\geq 0$ and $\sum_j t_j\leq 1$ then, for any $\varepsilon>0$ we have the next claim.

\smallskip

\noindent {\it Claim.} If we define the following sets:
\begin{align*}
A:= \left\{\omega \mid \, \exists \theta\in S^{k-1} :  \big| \| G_\omega(\theta)\|_{B_p(\mu)}^p -I_p^p \big| > \varepsilon I_p^p \right\}, \\
 A_1:= \left \{ \omega \mid \exists u_1\in {\cal N}_1 : \big| \|G_\omega(u_1)\|_{B_p(\mu)}^p -I_p^p \big| > t_1\varepsilon I_p^p \right\} \nonumber
 \end{align*} and for $j\geq 2$  
\begin{align*} 
A_j:= \left \{\omega \mid \exists u_j\in {\cal N}_j, u_{j-1}\in {\cal N}_{j-1} : \left| \|G_\omega(u_j)\|_{B_p(\mu)}^p - \|G_\omega(u_{j-1})\|_{B_p(\mu)}^p \right| > t_j \varepsilon I_p^p \right \},
\end{align*} where $I_p\equiv I_p(\gamma_n, B_p(\mu))$, then the inclusion
$A\subseteq \bigcup_{j=1}^\infty A_j$ holds.

\smallskip

\noindent {\it Proof of Claim.} If $\omega\notin \bigcup_{j=1}^\infty A_j$ then for any $j$ and any $u_j\in {\cal N}_j$ we have
\begin{align*}
\big | \|G_\omega(u_1)\|_{B_p(\mu)}^p - I_p^p \big| \leq \varepsilon t_1 I_p^p \quad {\rm and} \quad \left| \|G_\omega(u_j)\|_{B_p(\mu)}^p - \|G_\omega(u_{j-1})\|_{B_p(\mu)}^p \right| \leq \varepsilon t_j I_p^p, \quad j=2,3,\ldots .
\end{align*} For any $\theta$ there exist $u_j\in {\cal N}_j$ such that $\|\theta-u_j\|_2 <\delta_j$ for $j=1,2,\ldots$. Hence, for any 
$N\geq 2$ we may write
\begin{align*}
\big| \|G_\omega (\theta)\|_{B_p(\mu)}^p - I_p^p \big| &\leq \big| I_p^p-\|G_\omega(u_1)\|_{B_p(\mu)}^p \big| +\sum_{j=2}^N \left| \|G_\omega(u_{j-1})\|_{B_p(\mu)}^p-\|G_\omega(u_j)\|_{B_p(\mu)}^p \right| +\\
& + \big| \|G_\omega(u_N)\|_{B_p(\mu)}^p -\|G_\omega(\theta)\|_{B_p(\mu)}^p \big| \\
&\leq \sum_{j=1}^N \varepsilon t_j I_p^p +  p \cdot\delta_N\cdot \|G_\omega\|_{2\to X_p(\mu)}^p ,
\end{align*} which proves the assertion, since $N$ is arbitrary.

Fix $0 <\varepsilon<1$. Choose $\delta_j=e^{-j}$, $t_j = j^{p/2}e^{-j}/a_p$ with $a_p: =\sum_{j=1}^\infty j^{p/2}e^{-j}$ (thus, $\sum_j t_j \leq 1$). 
Then, if we employ the previous claim and Theorem \ref{thm:conc-B_p-2} we may write
\begin{align*}
\mathbb P(A) &\leq C|{\cal N}_1| \exp( -c_1 \psi(n,p,\varepsilon t_1) ) + C\sum_{j=2}^\infty |{\cal N}_{j-1}| \cdot |{\cal N}_j| \exp (- c_1\psi(n,p, \varepsilon t_je^j/4)) \\
&\leq C\sum_{j=1}^\infty (3 e^{j})^{2k} \exp \left( -c_1' \psi \left( n,p, \varepsilon j^{p/2} a_p^{-1} \right) \right) \\
& \leq C\sum_{j=1}^\infty \exp \left( c_2 j k- c_2' a_p^{-2} \psi \left( n,p,\varepsilon j^{p/2} \right)  \right)\\
& \leq C\sum_{j=1}^\infty \exp \left( c_2 j k - c_2' a_p^{-2} j \psi (n,p,\varepsilon) \right) \\
&\leq C \sum_{j=1}^\infty \exp \left( - c_3 a_p^{-2} j \psi(n,p,\varepsilon) \right) \leq C' \exp \left( - c_3' a_p^{-2} \psi(n,p,\varepsilon) \right)  ,
\end{align*} 
as long as $k\lesssim a_p^{-2} \psi(n,p,\varepsilon) \lesssim (2e/p)^p \psi(n,p, \varepsilon)$. Therefore, with probability greater than $1-e^{-ca_p^{-2} \psi(n,p,\varepsilon)}$ the random operator $G$ satisfies
\begin{align*}
(1-\varepsilon)^{1/p}I_p \|x\|_2 \leq \|G (x) \|_{B_p(\mu)} \leq (1+\varepsilon)^{1/p} I_p \|x\|_2,
\end{align*} for all $x \in \mathbb R^k$. To conclude we have to recall that ${\rm Im} G=F$ is Haar-distributed on $G_{n,k}$ 
(see \cite{Sch3} for the details). \prend

\smallskip

\begin{note} Let us mention that if $n^{-\frac{p-2}{2(p-1)}} \lesssim_p \varepsilon <1$ then we get $k\leq c (\varepsilon n)^{2/p} /p$
by taking into account the form of $\psi(n,p,\cdot)$. Indeed; for $n$ large enough, i.e. $n\geq e^{C p\log p}$, if we consider 
$(cp)^{p/2} n^{-\frac{p-2}{2(p-1)}} < \varepsilon <1$ in the previous series of inequalities we obtain
\begin{align*}
\mathbb P(A) \leq C\sum_{j=1}^\infty (3 e^{j})^{2k} \exp \left( -c_1' \psi \left( n,p, \varepsilon j^{p/2} a_p^{-1} \right) \right) 
&\leq C\sum_{j=1}^\infty \exp \left( c_2 j k- c_4 j (\varepsilon n)^{2/p} a_p^{-2/p}  \right) \\
&\leq C \sum_{j=1}^\infty \exp \left( c_2jk - c_5 j p^{-1} (\varepsilon n)^{2/p} \right) \\
&\leq C' \exp \left( -c_5' p^{-1} (\varepsilon n)^{2/p} \right),
\end{align*} provided that $k \leq c'p^{-1} (\varepsilon n)^{2/p} $.
\end{note}

%%%%%%%%%%%%%%%%%%%%%%%%%%%%%%%%%%%%%%%%%%%%%%%%%%%%%%%%%%%%
\section{Further remarks}
%%%%%%%%%%%%%%%%%%%%%%%%%%%%%%%%%%%%%%%%%%%%%%%%%%%%%%%%%%%%

\noindent {\bf 1. Optimality of the result.} If the isotropic measure $\mu$ on $S^{n-1}$ is the one supported on $\pm e_i$'s i.e.
$X_p(\mu)\equiv \ell_p^n$, then Theorem \ref{thm:conc-B_p} is optimal (up to constants depending on $p$) as was proved in \cite{PVZ}. 
Moreover, Theorem \ref{thm:main-3} is optimal, in the sense that if the 
typical $k$-dimensional subspace of $\ell_p^n$ is $(1+\varepsilon)$-spherical, then 
$k\leq C p (\varepsilon n)^{2/p} $ for some absolute constant $C>0$ (see \cite{PVZ}). We should mention that it is known, that for 
concrete values of $p$ one can embed $\ell_2^k$ into $\ell_p^n$ even isometrically (see \cite{Ko} for details). However, this
is not a typical subspace.

\bigskip

\noindent {\bf 2. Selection of randomness.} Embeddings of $\ell_2^k$ into $L_q, \; 2< q <\infty$ under different randomness have appeared in the 
literature in \cite{BDGJN}. The authors there consider large random matrices with independent Rademacher entries in order to $K(q)$-embed
$\ell_2^k$ into $\ell_q^N$ with $N\simeq k^{q/2}$, where $K(q)>0$ depends only on $q$. Then, they use this result in order 
to prove that for any $1<p<2$ there exists 
uncomplemented subspace of $L_p$ which is isomorphic to Hilbert space. It is worth mentioning, that one can prove a concentration 
result similar to that of Theorem \ref{thm:conc-B_p} using other randomness than Gaussian. 
In particular, if $\nu$ is an isotropic\footnote{The isotropic measures on $\mathbb R^n$ are defined similarly: A Borel probability measure $\nu$ on $\mathbb R^n$ is said to be isotropic if $\int_{\mathbb R^n} \langle x, \theta \rangle ^2 \, d\nu(x)=1$ for all $\theta\in S^{n-1}$.} Borel probability measure on $\mathbb R^n$ which satisfies a log-Sobolev inequality
with constant $\rho>0$ then we may prove the following:

\begin{theorem} \label{thm: gen-conc} Let $2<p<\infty$, let $\mu$ be a Borel isotropic measure on $S^{n-1}$ and let $\nu$ be an 
isotropic Borel probability measure on
$\mathbb R^n$ which satisfies a log-Sobolev inequality with constant $\rho>0$. Then, we have
\begin{align*}
\left( \iint \big| \|x\|_{B_p(\mu)}^p -\|y\|_{B_p(\mu)}^p \big|^r \, d\nu(x) d\nu(y) \right)^{1/r} \leq 
C(p, \rho) I_p^p(\nu, B_p(\mu)) \max \left\{ \left(\frac{r}{n}\right)^{1/2} , \, \frac{r^{p/2}}{n} \right\},
\end{align*} for all $r\geq 2$, where $C(p, \rho)>0$ is constant depending only on $p$ and $\rho$.
\end{theorem}

Having proved Theorem \ref{thm: gen-conc}, we apply Markov's inequality as in Section 3 to get the corresponding concentration inequality. 
For the proof of Theorem \ref{thm: gen-conc} we argue 
as follows: Consider the function $f(x)=\|x\|_{B_p(\mu)}^p$ and define $F=f-\mathbb E_\nu f$. Then, 
a direct application of Lemma \ref{lem:conc-log-sob} yields
\begin{align} \label{eq: 5.2}
\|F\|_{L_r(\nu)}^2 \leq \|F\|_{L_2(\nu)}^2 + \frac{1}{\rho} \int_2^r \big\| \|\nabla f\|_2 \big\|_{L_s(\nu)}^2 \, ds,
\end{align} for all $r\geq 2$. Recall the known fact (e.g. see \cite{Led}) that if a measure $\nu$ satisfies a log-Sobolev inequality with 
constant $\rho$, also satisfies a Poincar\'{e} inequality with constant $\rho$, that is
\begin{align*}
\|h -\mathbb E_\nu h\|_{L_2(\nu)}^2 \leq \frac{1}{\rho } \int_{\mathbb R^n} \| \nabla h\|_2^2 \, d\nu = \frac{1}{\rho} \big\| \|\nabla h\|_2 \big\|_{L_2(\nu)}^2,
\end{align*} for any smooth function $h$. Therefore, \eqref{eq: 5.2} becomes
\begin{align*}
\|F\|_{L_r(\nu)}^2 \leq \frac{2}{\rho} \int_2^r \big\| \|\nabla f\|_2 \big\|_{L_s(\nu)}^2 \, ds \leq \frac{2r}{\rho} \big\| \|\nabla f\|_2 \big\|_{L_r(\nu)}^2 ,
\end{align*} for all $r\geq 3$, where we have used the fact that $s \mapsto \|h\|_{L_s}$ is non-decreasing function. Taking into account 
the Claim in Section 3 we get
\begin{align} \label{eq: 5.5}
\|F\|_{L_r(\nu)}^2 \leq \frac{2p^2 r}{\rho}  \left( \int_{\mathbb R^n} \| x\|_{B_{2p-2}(\mu)}^{r(p-1)} \, d\nu(x) \right)^{2/r} , \quad r\geq 3.
\end{align} Again, Lemma \ref{lem:conc-log-sob} implies that
\begin{align*}
\left( \int_{\mathbb R^n} \| x\|_{B_{2p-2}(\mu)}^{r(p-1)} \, d\nu(x) \right)^{2/r} \leq
 I_{2p-2}^{2p-2}(\nu, B_{2p-2}(\mu))\left( 1+ \frac{(r-2)(p-1)}{\rho I_{2p-2}^2(\nu, B_{2p-2}(\mu))} \right)^{p-1},
\end{align*} for $r\geq 2$. Plug this back in \eqref{eq: 5.5} we obtain
\begin{align*}
\|F\|_{L_r(\nu)} < \left(\frac{2p^2}{\rho}\right)^{1/2} r^{1/2} I_{2p-2}^{p-1}(\nu, B_{2p-2}(\mu)) \left(1+\frac{(r-2)(p-1)}{\rho I_{2p-2}^2(\nu, B_{2p-2}(\mu))}\right)^{ \frac{p-1}{2}},
\end{align*} for all $r\geq 3$. Finally, we have
\begin{align*}
 n\leq I_p^p(\nu, B_p(\mu))  \leq \left(1+\frac{p-2}{\rho}\right)^{p/2} n,
\end{align*} for all $p\geq 2$. Indeed; we may write
\begin{align*}
I_p^p(\nu, B_p(\mu)) =\int_{S^{n-1}} \int_{\mathbb R^n} |\langle x, \theta \rangle|^p \, d\nu(x ) \, d\mu(\theta) 
\geq \int_{S^{n-1}} \left( \int_{\mathbb R^n} \langle x,\theta \rangle^2 \, \nu(x) \right)^{p/2}\, d\mu(\theta) = \mu(S^{n-1}),
\end{align*} where we have used H\"older's inequality and the isotropicity of $\nu$. 
For the right-hand side, we fix $\theta\in S^{n-1}$ and we apply Lemma \ref{lem:conc-log-sob} for $x\mapsto \langle x, \theta \rangle$
to get
\begin{align*}
\left ( \int_{\mathbb R^n} |\langle x,\theta \rangle|^p \, d\nu(x) \right)^{2/p} \leq 
\int_{\mathbb R^n} \langle x, \theta \rangle^2 \, d\nu(x)+\frac{p-2}{\rho}= 1+\frac{p-2}{\rho},
\end{align*} where we have used the isotropicity again. Finally, integration with respect to $\mu$ yields:
\begin{align*}
I_p^p(\nu, B_p(\mu)) =\int_{S^{n-1}} \int_{\mathbb R^n} |\langle x, \theta\rangle|^p \, d\nu(x) \, d\mu(\theta) \leq 
\int_{S^{n-1}} \left( 1+\frac{p-2}{\rho} \right)^{p/2} \, d\mu(\theta),
\end{align*} as asserted. Taking into account these estimates, we argue as in Section 3 to complete the proof. 
The details are left to the reader.

\bigskip

\noindent {\bf 3. Minimal Gaussian variance for subspaces of $L_p$.} Let us point out that our method also provides upper estimate for the variance 
of the norm of any 
finite dimensional subspace of $L_p$ in Lewis' position. We should mention that the following estimate turns out to be optimal
(up to constants of $p$) since they agree with the $\ell_p^n$ case (see \cite[Section 3]{PVZ} for details). 

\begin{theorem}
Let $1\leq p <\infty$. Then, for any $n$-dimensional subspace $X$ of $L_p$ represented on $\mathbb R^n$, 
there exists a position $\tilde B$ of its unit ball $B_X$ such that
\begin{align*}
{\rm Var}\|Z\|_{\tilde B} \leq C^p n^{\frac{2}{p}-1},
\end{align*} where $C>0$ is an absolute constant. In particular, the normalized variance is of minimal possible order (up to constants of $p$)
\begin{align*}
\frac{{\rm Var}\|Z\|_{\tilde B}}{\mathbb E\|Z\|^2_{\tilde B}} \leq \frac{C^p}{n}.
\end{align*}
\end{theorem}

\noindent {\it Sketch of Proof.} If $\tilde B$ is a Lewis' position of $B_X$, we may identify $X$ with $X_p(\mu)$ for some 
Borel isotropic measure $\mu$ on $S^{n-1}$. Then, the result follows from Corollary \ref{cor:var-B_p}. On the other hand note that
the normalized variance is minimal since for {\it every} norm $\|\cdot \|$ on $\mathbb R^n$ one has ${\rm Var}\|Z\| \gtrsim \mathbb E\|Z\|^2 /n$.
The latter may be easily checked by using integration in polar coordinates and the Cauchy-Schwarz inequality. \prend

\bigskip

\noindent {\bf 4. A Johnson-Lindenstrauss type result.} Note that using Lewis' lemma (Theorem \ref{thm:Lewis}) and the two level
concentration for the $B_p(\mu)$ bodies (Theorem \ref{thm:conc-B_p}) we can conclude the following low-dimensional embedding 
of Hilbertian sets to any subspace of $L_p,\; 2<p<\infty$. This can be viewed as a Johnson-Lindenstrauss type result for target spaces
which sit in $L_p$.

\begin{theorem} Let $2<p<\infty$. For any $\varepsilon\in (0,1)$ and for any $N\geq 1$, there exists 
$m\lesssim \max \left\{ p4^p\frac{\log N}{\varepsilon^2}, \frac{ (\log N)^{p/2}}{\varepsilon} \right\}$ with the following property: For any 
subset $S=\{x_1,\ldots,x_N\} \subset \ell_2$ and for any subspace $X$ of $L_p$ with $\dim X=m$, there exists a linear mapping 
$T: S\to X$ which satisfies
\begin{align*}
(1-\varepsilon) \|x_i-x_j\|_2 \leq \|Tx_i-Tx_j\|_X \leq (1+\varepsilon) \|x_i-x_j\|_2,
\end{align*} for all $i,j=1,\ldots,N$.
\end{theorem}

%\newpage
%%%%%%%%%%%%%%%%%%%%%%%%%%%%%%%%%%%%%%%%%%%%%%%%%%%%%%%%%%%%%%%%%%%%%%%%%

\bigskip

\vspace{.5cm} \noindent \begin{minipage}[l]{\linewidth}
  Grigoris Paouris: {\tt grigoris@math.tamu.edu}\\
  Department of Mathematics, Mailstop 3368\\
  Texas A\&M University\\
  College Station, TX 77843-3368\\
  
  \bigskip
  
  Petros Valettas: {\tt valettasp@missouri.edu}\\
  Mathematics Department\\
  University of Missouri\\ 
  Columbia, MO 65211\\

\end{minipage}

\end{document}